\documentclass[journal]{IEEEtran}
\let\labelindent\relax
\usepackage{enumitem}
\setlist[itemize]{leftmargin=*}
\usepackage{graphicx}
\usepackage{array}
\usepackage{textcomp}
\usepackage{cite}
\usepackage{amsmath}
\usepackage{multirow}
\usepackage{bm}
\usepackage{verbatim}
\usepackage{xcolor}
\usepackage{threeparttable}
\usepackage{url}

\ifCLASSOPTIONcompsoc
\usepackage[caption=false,font=normalsize,labelfont=sf,textfont=sf]{subfig}
\else
\usepackage[caption=false,font=footnotesize]{subfig}
\fi

\ifCLASSINFOpdf
\else
\fi
\begin{document}
%
\title{Security Constrained Multi-Stage Transmission Expansion Planning Considering a Continuously Variable Series Reactor}
%
%
%

\author{Xiaohu~Zhang,~\IEEEmembership{Student Member,~IEEE,}
        Kevin~Tomsovic,~\IEEEmembership{Fellow,~IEEE,} \\
        and Aleksandar~Dimitrovski,~\IEEEmembership{Senior Member,~IEEE}
\thanks{This work was supported in part by ARPAe (Advanced Research Projects
	Agency Energy), and in part by the Engineering Research Center Program
	of the National Science Foundation and the Department of Energy under
	NSF Award Number EEC-1041877 and the CURENT Industry Partnership Program.}
\thanks{Xiaohu Zhang and Kevin Tomsovic are with the Department of Electrical Engineering and Computer Science, the University of Tennessee, Knoxville, TN, USA, email:\{xzhang46,tomsovic\}@utk.edu. Aleksandar Dimitrovski is with the Department of Electrical and Computer Engineering, University of Central Florida, FL, USA, email: Aleksandar.Dimitrovski@ucf.edu.}}

%
%

\markboth{PREPRINT OF DOI: 10.1109/TPWRS.2017.2671786, IEEE TRANSACTIONS ON POWER SYSTEMS}%
{Shell \MakeLowercase{\textit{et al.}}: Bare Demo of IEEEtran.cls for Journals}
%



\maketitle

\begin{abstract}
This paper introduces a Continuously Variable Series Reactor (CVSR) to the transmission expansion planning (TEP) problem. The CVSR is a FACTS-like device which has the capability of controlling the overall impedance of the transmission line. However, the cost of the CVSR is about one tenth of a similar rated FACTS device which potentially allows large numbers of devices to be installed. The multi-stage TEP with the CVSR considering the $N-1$ security constraints is formulated as a mixed integer linear programming model. The nonlinear part of the power flow introduced by the variable reactance is linearized by a reformulation technique. To reduce the computational burden for a practical large scale system, a decomposition approach is proposed. The detailed simulation results on the IEEE 24-bus and a more practical Polish 2383-bus system demonstrate the effectiveness of the approach. Moreover, the appropriately allocated CVSRs add flexibility to the TEP problem and allow reduced planning costs. Although the proposed decomposition approach cannot guarantee global optimality, a high level picture of how the network can be planned reliably and economically considering CVSR is achieved.    
\end{abstract}


%
\IEEEpeerreviewmaketitle


\section*{Nomenclature}
\subsection*{Indices}
\addcontentsline{toc}{section}{Nomenclature}
\begin{IEEEdescription}[\IEEEusemathlabelsep\IEEEsetlabelwidth{$V_1,V_2,V_3$}]
\item[$i, \ j$] Index of buses.
\item[$k$] Index of transmission elements.
\item[$n$] Index of generators.
\item[$m$] Index of loads.
\item[$c$] Index of states; $c=0$ indicates the base case; $c>0$ is a contingency state.
\item[$b$] Index of load block
\item[$t$] Index of time.
\item[$E$] Index for an existing transmission line.
\item[$C$] Index for a candidate transmission line.
\end{IEEEdescription}

\subsection*{Variables}
\addcontentsline{toc}{section}{Nomenclature}
\begin{IEEEdescription}[\IEEEusemathlabelsep\IEEEsetlabelwidth{$V_1,V_2,V_3$}]
\item[$P^g_{ncbt}$] Active power generation of generator $n$ for state $c$ in load block $b$ at time $t$.
\item[$P_{kcbt}$] Active power flow on branch $k$ for state $c$ in load block $b$ at time $t$.
\item[$x^V_{k}$] Reactance of a CVSR at branch $k$.
\item[$u_{k,1}^E,u_{k,2}^E$] Slack variables for the flow violation at existing branch $k$.
\item[$u_{k,1}^C,u_{k,2}^C$] Slack variables for the flow violation at candidate branch $k$.
\item[$\theta_{kcbt}$] The voltage angle difference across the branch $k$ in load block $b$ for state $c$ at time $t$.
\item[$\alpha_{kt}$] Binary variable associated with line investment for branch $k$ at time $t$.
\item[$\delta_{kt}$] Binary variable associated with placing a CVSR on branch $k$ at time $t$.
\end{IEEEdescription} 
\subsection*{Parameters}
\addcontentsline{toc}{section}{Nomenclature}
\begin{IEEEdescription}[\IEEEusemathlabelsep\IEEEsetlabelwidth{$V_1,V_2,V_3$}]
\item[$b_k$] Susceptance for branch $k$.
\item[$d$] Discount factor.
\item[$n_l$] Number of branches in the system.
\item[$x_k$] Reactance for branch $k$.
\item[$x_{k,V}^{\min}$] Minimum reactance of the CVSR at branch  $k$.
\item[$x_{k,V}^{\max}$] Maximum reactance of the CVSR at branch  $k$.
\item[$A_{bt}$] Duration of the load block $b$ at time $t$.
\item[$C_k^{V}$] Investment cost of the CVSR at branch $k$.
\item[$C_k^L$] Investment cost of the branch $k$.
\item[$C_n^g$] Operating cost coefficient for generator $n$.
\item[$N_{kcbt}$] Binary parameter associated with the status of branch $k$ for state $c$ in load block $b$ at time $t$. 
\item[$P_{ncbt}^{g,\min}$] Minimum active power output of generator $n$ for state $c$ in load block $b$ at time $t$.
\item[$P_{ncbt}^{g,\max}$] Maximum active power output of generator $n$ for state $c$ in load block $b$ at time $t$.
\item[$P_{mcbt}^d$]  Active power consumption of demand $m$ for state $c$ in load block $b$ at time $t$.
\item[$S_{kcbt}^{\max}$] Thermal limit of branch  $k$ for state $c$ in load block $b$ at time $t$. 
\item[$\theta_k^{\max}$] Maximum angle difference across branch $k$: $\pi/3$ radians.

\end{IEEEdescription} 

\subsection*{Sets}
\addcontentsline{toc}{section}{Nomenclature}
\begin{IEEEdescription}[\IEEEusemathlabelsep\IEEEsetlabelwidth{$V_1,V_2,V_3$}]
\item[$\mathcal{D}_i$] Set of loads located at bus $i$.
\item[$\Omega_{L}$] Set of existing transmission lines.
\item[$\Omega_{L}^+$] Set of candidate transmission lines.
\item[$\Omega_{L}^i$] Set of transmission lines connected to bus $i$.
\item[$\Omega_{T}$] Set of time periods.
\item[$\Omega_{c}$] Set of states.
\item[$\Omega_{b}$] Set of load blocks.
\item[$\Omega_{0}$] Set of base operating state.
\item[$\Omega_{V}$] Set of candidate transmission lines to install CVSR.
\item[$\mathcal{B}$] Set of buses. 
\item[$\mathcal{G}$] Set of on-line generators.
\item[$\mathcal{G}_{i}$] Set of on-line generators located at bus $i$.
\item[$\mathcal{G}_{fix}$] Set of on-line generators with fixed generation. 
\end{IEEEdescription}

\section{Introduction}
\IEEEPARstart{T}{he} Continuously Variable Series Reactor (CVSR) has recently been proposed for power flow control \cite{mybibb:Aleks,mybibb:SCR}. By controlling the saturation of a magnetic core, the device is capable of continuously and smoothly regulating its output reactance, which is similar to a series FACTS controller TCSC. The control circuit for the CVSR is a simple and low power rating AC/DC converter so the cost of the CVSR is far less than that of the TCSC. Numerous CVSRs could be installed into a single system to enable comprehensive use of the transmission capacity. This could have a significant impact on Transmission Expansion Planning (TEP) decisions and is the main reason for revisiting the TEP problem formulation in this paper.

TEP is a task that determines the best strategy to add new transmission lines to the existing power network in order to satisfy the growth of electricity demand and generation over a specified planning horizon. In the contemporary power system, due to the power market restructuring and massive integration of renewable energy, it is critical to have a rationally planned power system that is not only capable of serving the increasing load reliably and efficiently but also economically \cite{mybibb:TEP_huizhang1}.  Depending on the model, TEP can be classified as either a single-stage or multi-stage model. For a single-stage TEP, additional lines are planned only for the target planning year; while for the multi-stage TEP, several different planning horizons with distinct load and generation patterns are considered together. Multi-stage TEP not only decides \textit{where} to build the new transmission line, but also determines \textit{when} to build the new line \cite{mybibb:Jabr_TEP,mybibb:TEP_strategy}.

The modeling and solution techniques for the traditional TEP problem have been studied extensively. Mathematical programming is a major category of the solution methods. At the transmission level, the DC power flow model is capable of providing a good approximation and linear methods can be applied. In \cite{mybibb:TEP_2003,mybibb:TEP_letter}, the TEP in DC network model was formulated as a mixed integer linear programming (MILP) problem and solved by a commercial optimization solver. A disjunctive factor was introduced to eliminate the product between continuous and  binary variables. Given the non-convex nature of the power system, the exact AC network model for the TEP problem is generally a non-convex mixed integer nonlinear programming (MINLP) problem. This type of model is challenging for existing commercial solvers. Therefore, several relaxed or approximated AC models for the TEP problem have been proposed. 

In \cite{mybibb:TEP_huizhang2,mybibb:TEP_Taylor}, the nonlinear AC power flow equations were linearized around the operating point based on Taylor series to achieve the linear model for the AC TEP. The quadratic constraints, such as, the active and reactive power losses, the MVA limit for the transmission line were approximated by using piecewise linearization. In \cite{mybibb:TEP_Josua1}, the lift and project \cite{mybibb:lift_project1} technique was adopted to lift the TEP problem into higher dimensional space and project the relaxed solution onto the original space. In \cite{mybibb:TEP_Joshua}, the line flow based power flow equations \cite{mybibb:LFB} were employed to give a convex second order cone model for the AC TEP. The voltage magnitude was assumed to be equal to one and the non-convex constraint for the voltage drop across a transmission line was omitted. The AC or relaxed AC TEP models provide a relatively more accurate representation of the network and can include the reactive power planning (RPP) into the TEP problem. However, to the best of the authors' knowledge, the AC TEP models were only applied to small or medium scale systems. Meta-heuristic methods, such as, genetic algorithms \cite{mybibb:TEP_GA}, greedy randomized search \cite{mybibb:TEP_greedy_random}, particle swarm optimization \cite{mybibb:TEP_PSO} and differential evolution \cite{mybibb:TEP_DE} have also been proposed to solve the TEP problem. These techniques have the advantage of easy and straightforward implementation; however, they suffer disadvantages of susceptibility to local optimum and slow computational speed for large practical systems \cite{mybibb:heuristic_review1,mybibb:heuristic_review2}. 

Major hurdles for construction of new transmission lines are difficulties in obtaining the right-of-way, political resistance, long construction time and limited capital budget. These challenging issues have drawn interest in techniques for delaying upgrades. In \cite{mybibb:TS_TEP}, transmission switching (TS) was introduced to defer the construction of new transmission lines. Benders Decomposition was employed to solve the planning and operation problem alternately. In \cite{mybibb:TEP_FACTS}, the authors evaluated the economic benefits and increased flexibility by including the FACTS devices in the TEP. In \cite{mybibb:TEP_ES}, a single stage TEP model considering energy storage systems (ESS) was presented. The total investment cost for the transmission lines can be reduced by appropriately placing the ESS in the system. 

This paper presents a MILP model for the multi-stage TEP considering CVSRs, while satisfying $N-1$ security constraints. Three load blocks are selected to accommodate the load profile of each stage and the considered transmission contingencies can occur in any of the load blocks. Several benefits are anticipated by introducing the CVSR into TEP: 1) CVSRs improve the utilization of the existing network, which leads to deferment or avoidance of new transmission lines; 2) CVSRs change the power flow pattern and increase the use of lower cost  generation, which reduces the total operating cost; 3) CVSRs add flexibility to the system and provide additional corrective actions following contingencies.  

The main contributions of this paper are summarized below:
\begin{enumerate}
	\item A security constrained multi-stage TEP with the consideration of CVSRs is formulated.  
	\item A reformulation technique is proposed to transform the MINLP model into a MILP model allowing the problem to be solved by mature commercial MILP solvers.
	\item An iterative approach is developed to decompose the model into the planning master problem and the security check sub-problem so that it is computationally tractable for practical sized systems. This is critical as the model size increases dramatically with the number of stages, load blocks and contingencies. 
\end{enumerate} 
 
Due to the heuristic method used in the iterative approach, the solution obtained by the decomposition model is not guaranteed to be global optimality. However, it provides a high level picture of how the network can be rationally planned including CVSRs so it is useful from engineering point of view. In addition, the decomposition approach allows the originally large scale MINLP model to become tractable.

The remainder of the paper is organized as follows. In Section \ref{static_model}, the steady state model of CVSR in DC power flow is presented and the reformulation technique is illustrated to transform the originally nonlinear power flow model into a linear model. Section \ref{Optimization_model} presents detailed information about the optimization model and the iterative approach. Simulation results are given in Section \ref{case_study} on the IEEE 24-bus and a more practical Polish systems. Conclusions are given in Section \ref{conclusions}.

\section{Steady State Model of CVSR and the Reformulation Technique}
\label{static_model}
\subsection{Steady State Model of CVSR in DC Power Flow}
\begin{figure}[!htb]
\centering
\includegraphics[width=0.4\textwidth]{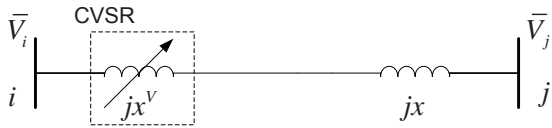}
\caption{Static representation of CVSR in DC power flow.}
\label{CVSR_steady}
\end{figure}

Fig. \ref{CVSR_steady} depicts the usage of a CVSR as a series control device. It can be represented by a continuously variable inductive reactance with the parasitic resistance ignored.

With a CVSR inserted on a transmission line and the resistance ignored, the total susceptance of the transmission line becomes:
\begin{equation}
b_k^{\prime}=-\frac{1}{x_k+x_k^{V}}=-(b_{k}+b_k^{V})
\end{equation}
where 
\begin{align}
&b_{k}=\frac{1}{x_k}   \\
&b_k^{V}=-\frac{x_k^{V}}{x_k(x_k+x_k^{V})}
\end{align}

Unlike a TCSC, the CVSR can only provide a positive reactance and here it is assumed that the CVSR will be installed on the transmission line which is not overly compensated by a series capacitor, i.e., $x_k \ge 0$. Then the range of the variable susceptance $b_{k}^V$ is:
\begin{align}
b_{k,V}^{\min}&=-\frac{x_{k,V}^{\max}}{x_k(x_k+x_{k,V}^{\max})}  \label{bCVSRmin} \\
b_{k,V}^{\max}&=-\frac{x_{k,V}^{\min}}{x_k(x_k+x_{k,V}^{\min})}  \label{bCVSRmax}
\end{align}

\subsection{Reformulation of the Nonlinear Power Flow Equation}
The active power flow on the candidate transmission line assuming a DC power flow model can be expressed as:
\begin{align}
&P_k=(b_{k}+\delta_kb_k^{V})\cdot \theta_k  \label{power_CVSR}  \\
&b_{k,V}^{\min} \le b_k^{V} \le b_{k,V}^{\max}   \label{bcvsr_range} 
\end{align}

In (\ref{power_CVSR}), the binary variable $\delta_k=1$ indicates that a CVSR is installed on line $k$. The nonlinearity in (\ref{power_CVSR}) results from the trilinear term $\delta_kb_k^{V}\theta_k$. To linearize the nonlinear part, a new variable $w_k$ is introduced as:
\begin{equation}
w_k=\delta_kb_k^{V}\theta_k  \label{wij1}
\end{equation}

Then the active power flow equation (\ref{power_CVSR}) can be given as:
\begin{equation}
P_k=b_{k}\theta_k+w_k  \label{PCVSR_wij1}
\end{equation}

We multiply each side of the constraint (\ref{bcvsr_range}) with the binary variable $\delta_k$ and combine with the variable $w_k$:
\begin{equation}
\delta_k b_{k,V}^{\min} \le \frac{w_k}{\theta_k}=\delta_kb_k^{V} \le \delta_k b_{k,V}^{\max}  \label{if_ineq}
\end{equation}

The sign of $\theta_k$ determines the allowable range for $w_k$:
\begin{equation}
\left\lbrace 
\begin{aligned}
\delta_k\theta_k b_{k,V}^{\min} \le w_k \le \delta_k\theta_k b_{k,V}^{\max}, \ &\text{if} \ \theta_k>0  \\
w_k=0,\ \ \ \ \ \ \ \ \ \ \ \ \ \ \ \ \ \ \ \ \ \ \ \ \ &\text{if} \ \theta_k=0   \\
\delta_k\theta_k b_{k,V}^{\max} \le w_k \le \delta_k\theta_k b_{k,V}^{\min}, \ &\text{if} \ \theta_k<0 
\end{aligned} \right.
\end{equation}

To realize the ``if" constraints, an additional binary variable $y_k$ and the big-M complementary constraints \cite{mybibb:Tao_bigM} are introduced:
\begin{equation}
-M_ky_k+\delta_k\theta_kb_{k,V}^{\min} \le w_k 
\le \delta_k\theta_kb_{k,V}^{\max}+M_ky_k \label{if_1} \\
\end{equation}
\begin{align}
-M_k(1-y_k)+\delta_k\theta_k&b_{k,V}^{\max} \le w_k \nonumber  \\ 
 &\le \delta_k\theta_kb_{k,V}^{\min}+M_k(1-y_k) \label{if_2} 
\end{align}

During the optimization process, only one of the two constraints (\ref{if_1}) and (\ref{if_2}) will become active and the other one will be a redundant constraint that is always satisfied with a sufficiently large number $M_k$. Specifically, when $\theta_k< 0$, $y_k$ will be equal to one and the constraint (\ref{if_2}) will be active; when $\theta_k > 0$, $y_k$ will be equal to zero and the constraint (\ref{if_1}) will be active; when $\theta_k=0$, one of these two constraints will drive $w_k$ to zero regardless of the value of $y_k$. Note that the numerical problems occur when $M_k$ is chosen to be too large \cite{mybibb:TS_contingency,mybibb:investment_naps}. Since $b_{k}^V$ is negative, $M_k$ is selected as $|b_{k,V}^{\min}\theta_k^{\max}|$.

Now in the constraints (\ref{if_1}) and (\ref{if_2}), it can be seen that there is still a nonlinear term $\delta_k\theta_k$ which involves the product of a binary variable and a continuous variable. We introduce another new variable $z_k=\delta_k\theta_k$ and linearize it using the method in \cite{mybibb:WP_Adams}:
\begin{align}
&-\delta_k\theta_k^{\max} \le z_k \le \delta_k\theta_k^{\max} \label{z1} \\
&\theta_k-(1-\delta_k)\theta_k^{\max} \le z_k \le \theta_k+(1-\delta_k)\theta_k^{\max} \label{z2}
\end{align}

We then substitute $\delta_k\theta_k$ with $z_k$ in the inequalities (\ref{if_1}) and (\ref{if_2}):
\begin{align}
&-M_ky_k+z_kb_{k,V}^{\min} \le w_k 
\le z_kb_{k,V}^{\max}+M_ky_k \label{if3} \\
&-M_k(1-y_k)+z_kb_{k,V}^{\max} \le w_k\le z_kb_{k,V}^{\min}+M_k(1-y_k)\label{if4} 
\end{align}

Thus, the power flow equations (\ref{power_CVSR}) and (\ref{bcvsr_range}) in MINLP model is transformed to a MILP model with (\ref{PCVSR_wij1}) and (\ref{z1})-(\ref{if4}).

\section{Optimization Model}
\label{Optimization_model}
\subsection{$N-1$ Security Constraints}
Power grid security is the primary concern for the system operations and planning and it cannot be compromised. According to the NERC planning standards \cite{mybibb:NERC_standard}, a rationally planned power system should have the capability of maintaining an $N-1$ secure network. To include the $N-1$ contingency for the transmission lines into the optimization model, a binary parameter $N_{kc}$, which represents the status of line $k$ in state $c$ is introduced \cite{mybibb:TS_AC}:
\begin{equation}
N_{kc}=\left\lbrace 
\begin{aligned}
& 1, \ \text{if line}\ k\ \text{is in service in state}\ c  \\
& 0, \ \text{if line}\ k\ \text{is out of service in state}\ c 
\end{aligned} \right.
\end{equation}
It should be noted that $N_{k0}$ is equal to one since no transmission element is in outage for the base operating condition. The number of states for a complete transmission line $N-1$ contingency and the base case is $n_l+1$.  

For most planning problems, a complete set of $N-1$ contingencies is not needed and just results in excessive computations as $n_l$ is large in a practical system. For the TEP problem, a complete $N-1$ contingency is not needed since the addition of some new transmission lines in one area will mainly affect the power flow pattern in the nearby areas. The selection of the contingencies can be based on experimental data or a contingency screening algorithms \cite{mybibb:Cont_screening,mybibb:TEP_huizhang1}.

\subsection{Integrated Planning Formulation}
\label{integrated_formulation}
The integrated planning indicates that all the planning stages, load blocks and security constraints are included in one planning problem, which is formulated as (\ref{obj})-(\ref{time_couple2}).

\subsubsection{Objective Function}
The objective employed in this paper for the TEP problem minimizes the total cost, which includes both the investment and operating cost. Assuming a fixed load demand (price inelastic), minimizing operating cost is equivalent to minimizing generation cost. The objective function is:
\begin{align}
\min \ &\sum_{t \in TPL}\sum_{k \in \Omega_{L}^{+}}\frac{C_k^{L}(\alpha_{kt}-\alpha_{k,t-1})}{(1+d)^{t-1}} \nonumber \\
& +\sum_{t \in TPL}\sum_{k \in \Omega_{V}}\frac{C_k^{V}(\delta_{kt}-\delta_{k,t-1})}{(1+d)^{t-1}} \nonumber \\
& +\sum_{t \in TPL}\sum_{b \in \Omega_b}\sum_{n \in \mathcal{G}}\frac{A_{bt}C^g_{n}P^g_{n0bt}}{(1+d)^{t-1}}  \label{obj}
\end{align}

TPL represents the total planning horizon. The first two terms represent the one time investment cost for the new transmission lines and the installed CVSRs. The third term is the generation cost across the operating horizon. Three distinct load patterns which represent peak, normal and low load condition are selected to accommodate the load profile in each stage. Here the generation cost is just an estimated cost. However, if the detailed load duration curve for each year is given, a relatively more accurate generation cost model can be formulated. All the cost terms are discounted to the present value by using the discount factor $d$. In this paper, $d$ is selected to be 5\%.

\subsubsection{Constraints}
The active power flow through the existing transmission lines is:
\begin{align}
& P_{kcbt}^E-b_{k}\theta_{kcbt}+M_k'(1-N_{kcbt}) \ge 0,k\in \Omega_{L}\backslash \Omega_{V} \label{norm_c11}  \\
& P_{kcbt}^E-b_{k}\theta_{kcbt}-M_k'(1-N_{kcbt}) \le 0,k\in \Omega_{L}\backslash \Omega_{V} \label{norm_c22}  \\
& P_{kcbt}^E-b_{k}\theta_{kcbt}-w_{kcbt}+M_k'(1-N_{kcbt}) \ge 0, k\in \Omega_{V} \label{CVSR_c11} \\
& P_{kcbt}^E-b_{k}\theta_{kcbt}-w_{kcbt}-M_k'(1-N_{kcbt}) \le 0, k\in \Omega_{V} \label{CVSR_c22} 
\end{align}
Constraints (\ref{norm_c11})-(\ref{CVSR_c22}) hold $\forall c \in \Omega_c, b \in \Omega_b, t \in \Omega_T $.

Constraints (\ref{norm_c11}) and (\ref{norm_c22}) denote the active power on the lines without CVSRs while constraints (\ref{CVSR_c11}) and (\ref{CVSR_c22}) represent the active power flow on the candidate lines to install CVSRs. If the line is in service, i.e. $N_{kcbt}=1$, the line flow equations are enforced. A large disjunctive factor $M_k'$ is introduced to ensure these constraints are not restrictive when the transmission line is out of service. As the phase angle will not fall outside of the range $[-\pi/2 \ \pi/2 ]$ if an appropriate slack bus is selected, $M_k'$ is chosen to be $|b_{k}\pi|$.

Additional constraints introduced by the reformulation technique can be expanded to consider multiple states, load blocks and stages:
\begin{align}
&-M_ky_{kcbt}+z_{kcbt}b_{k,V}^{\min} \le w_{kcbt} 
\le z_{kcbt}b_{k,V}^{\max}+M_ky_{kcbt} \label{big_M11} \\
&-M_k(1-y_{kcbt})+z_{kcbt}b_{k,V}^{\max} \le w_{kcbt} \nonumber \\
&\ \ \ \ \ \ \ \ \ \ \ \ \ \ \ \ \ \ \ \ \ \le z_{kcbt}b_{k,V}^{\min}+M_k(1-y_{kcbt}) \label{big_M22}  \\
& -N_{kcbt}\delta_{kt}\theta_k^{\max} \le z_{kcbt} \le N_{kcbt}\delta_{kt}\theta_k^{\max} \label{bilinear_c11} \\
& N_{kcbt}(\theta_{kcbt}-(1-\delta_{kt})\theta_k^{\max}) \le z_{kcbt}   \nonumber   \\
&\ \ \ \ \ \ \ \ \ \ \ \ \ \ \ \ \ \ \ \ \le N_{kcbt}(\theta_{kcbt}+(1-\delta_{kt})\theta_k^{\max}) \label{bilinear_c22}
\end{align}
Constraints (\ref{big_M11})-(\ref{bilinear_c22}) hold $\forall k\in \Omega_{V},c \in \Omega_c,b \in \Omega_b, t \in \Omega_T$.

Constraints (\ref{big_M11})-(\ref{bilinear_c22}) guarantee that the line flow change $w_{kcbt}$ introduced by the CVSR is zero when  line $k$ with CVSR is out of service in state $c$, load block $b$ and at stage $t$. 

The power flow through the candidate transmission lines is:
\begin{align}
& P_{kcbt}^C-b_{k}\theta_{kcbt}+M_k'(2-N_{kcbt}-\alpha_{kt}) \ge 0 \label{can_c11}  \\
& P_{kcbt}^C-b_{k}\theta_{kcbt}-M_k'(2-N_{kcbt}-\alpha_{kt}) \le 0 \label{can_c22} 
\end{align}
Constraints (\ref{can_c11})-(\ref{can_c22}) hold $\forall k\in \Omega_{L}^+,c \in \Omega_c,b \in \Omega_b, t \in \Omega_T$.

In contrast with the existing transmission lines, a candidate transmission line has two situations where it is not connected: either it is not built or it has been built but is out of service.

The active power nodal balance at each bus is:
\begin{align}
&\sum_{n \in \mathcal{G}_i}P_{ncbt}^g -\sum_{m \in \mathcal{D}_i}P_{mcbt}^d=\sum_{k\in \Omega_{L}^{i} } P_{kcbt}^E+\sum_{k\in \Omega_L^{i} } P_{kcbt}^C  \label{P_B1} \\ 
& \ \ \ \ \ \ \ \ \ \ \ \ \ \ \ \ \ \ \ \ \ \  \ \ \ \ \ \ \ \ \ i \in \mathcal{B}, c \in \Omega_c,b \in \Omega_b, t \in \Omega_T \nonumber
\end{align}

The system physical limits are represented by:
\begin{align}
& -N_{kcbt}S_{kcbt}^{\max} \le P_{kcbt}^E \le N_{kcbt} S_{kcbt}^{\max},k \in \Omega_L \label{Slimit_E}  \\ 
& -\alpha_{kt} N_{kcbt}S_{kcbt}^{\max} \le P_{kcbt}^C \le \alpha_{kt}N_{kcbt} S_{kcbt}^{\max}, k \in \Omega_L^+ \label{Slimit_C} \\
& P_{ncbt}^{g,\min}\le P_{ncbt}^g \le P_{ncbt}^{g,\max}, n \in \mathcal{G}   \label{Pg_limit1} \\
&P_{ncbt}^g=P_{n0bt}^g, \  n \in \mathcal{G}_{fix}, c \in \Omega_c \backslash \Omega_{0}, b \in \Omega_b,t \in \Omega_T  \label{Pg_fix} 
\end{align}
Constraints (\ref{Slimit_E})-(\ref{Pg_limit1}) hold $\forall c \in \Omega_c, b \in \Omega_b,t \in \Omega_T$. Constraints (\ref{Slimit_E}) and (\ref{Slimit_C}) ensure that the power flow is zero if the line is not built or out of service; otherwise, the power flow on the line is limited by its thermal rating. Constraints (\ref{Pg_limit1}) and (\ref{Pg_fix}) reflect that only a subset of the generators are allowed to re-dispatch after a contingency. The other generators which do not participate in the rescheduling are fixed at their base case power output.

The build decisions made in the current stage must be present on the later stage:
\begin{align}
& \alpha_{kt} \ge \alpha_{k,t-1}, \ k \in \Omega_L^+,t \in \Omega_T  \label{time_couple1} \\
& \delta_{kt} \ge \delta_{k,t-1}, \ k \in \Omega_{V},t \in \Omega_T  \label{time_couple2} 
\end{align}
Note that $\alpha_{k0}$ and $\delta_{k0}$ are set to be zero.


\subsection{Decomposition}
\begin{figure}[!htb]
	\centering
	\includegraphics[width=0.385\textwidth]{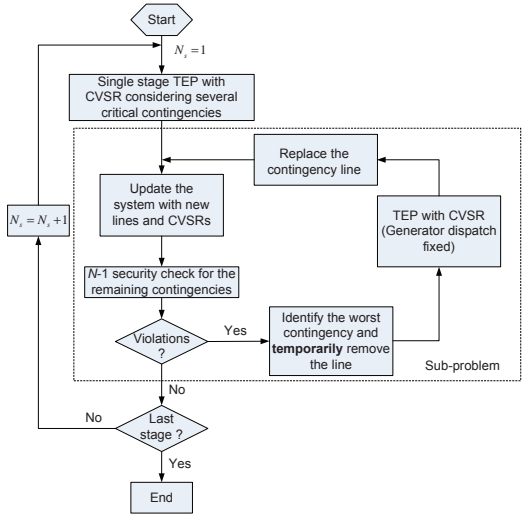}
	\caption{Flow chart of the iterative approach.}
	\label{decomposed}
\end{figure}

In the integrated planning model, the constraints have four dimensions, i.e., power system element, state, load block and time. Hence, the size of the optimization model will grow dramatically with the system size and planning horizon. To reduce the  computational burden for a large practical planning problem, the multiple stages are decomposed using forward planning \cite{mybibb:TEP_letter,mybibb:Jabr_TEP}, in which the planning for each stage is solved successively while the building decisions from the previous stage are enforced on subsequent stages. Although forward planning may lead to a suboptimal plan, it greatly reduces the computational time with relatively minor degradation of the solution quality. This iterative approach is as depicted in Fig. \ref{decomposed}.

 Essentially the majority of the $N-1$ security analysis will be performed iteratively at the sub-problem level. Checking $N-1$ security constraints iteratively is an effective way to reduce the computational burden in the TEP problem \cite{mybibb:TEP_huizhang1,mybibb:cc_iterative1}. The majority of the utilities use similar approaches for solving security constrained TEP \cite{mybibb:cc_iterative2}. The process is as below:
\begin{enumerate}
	\item Initialization of the stage number $N_s=1$. 
	\item Run the single stage TEP with CVSR model for the base case considering all the load blocks and several critical contingencies (CC). Obtain solutions and update the system with the new transmission lines and CVSRs.  
	\item Perform the remaining $N-1$ security analysis for the expanded system. If there are no violations, go to step 5); otherwise, identify the contingency leading to the worst violations. \textbf{Temporarily} remove the line from the system.
	\item Run the TEP with CVSR model. The generation dispatch is assumed to be unchanged. The purpose of this step is to find the optimal building plan (lines and CVSRs) to resolve the worst contingency. Replace the contingency line and update the system with new lines and CVSRs from this solution, go to step 3).
	\item If the last stage is solved, then complete; otherwise, increase the stage number $N_s=N_s+1$ and go to step 2). 
\end{enumerate}

Including several critical contingencies in the master problem is motivated by the natural thought that the critical contingencies have large impacts on the TEP results. However, considering more contingencies tends to increase the dimension of the master problem. The computational issues are discussed in Section \ref{computational}. For a practical system, the critical contingencies can be selected based on  empirical data. In our test system, we rank the contingencies in terms of circuit loading and the generation cost \cite{mybibb:automatic_contingency}. 

The two sections below detail the problem formulation of the master problem and sub-problem described above. Note that the constraints in Section \ref{integrated_formulation} all pertain to a specific state $c$, load block $b$ and stage $t$ . 
 
\subsubsection{Master Problem}
The planning master problem is to obtain the optimal building plan for the base case considering several critical contingencies. The optimization minimizes (\ref{obj}) subject to (\ref{norm_c11})-(\ref{time_couple2}). Note that the solution from the previous stage is the input for the current stage, i.e., $\alpha_{k,t-1}$ and $\delta_{k,t-1}$ are known before solving stage $t$.

\subsubsection{Sub-problem}
After obtaining the solution for the master problem in stage $t$, the sub-problem performs $N-1$ security analysis for the expanded system. Here, $P_{n0bt}^g$, $\alpha_{kt}$ and $\delta_{kt}$ are all input values for the security sub-problem while $P_{n0bt}^g$ is from the base case generation for each load block. In the iterative process of the sub-problem, new lines and CVSRs will be added to resolve the contingency, i.e., step 4), so $\alpha_{kt}$ and $\delta_{kt}$ need to be updated accordingly at each iteration.  

The violations for the DC power flow model are only thermal limit violations. For the $N-1$ security check, we introduce four positive slack variables to represent possible violations of the existing and candidate transmission lines. For each contingency state $c$, the objective is to minimize the sum of these slack variables:
\begin{equation}
\min \ \sum_{k\in \Omega_L}(u_{k,1}^E+u_{k,2}^E)+\sum_{k\in \Omega_L^+}(u_{k,1}^C+u_{k,2}^C)
\end{equation}
Obviously, the contingency with the maximum objective will be regarded as the worst contingency. If there is no violation, the objective for all the contingencies must fall within a specified tolerance.
The thermal limit constraints are:
\begin{align}
& -N_k(S_k^{\max}+u_{k,1}^E) \le P_k^E \le N_k(S_k^{\max}+u_{k,2}^E), \nonumber \\
&\ \ \ \  \ \ \ \ \ \ \ \ \ \ \ \ \ \ \ \ \ \ \ \ \ \ \ \ \ \ \ \ \ \ \ \ \ \ \ \ \ \ \ \ \ \ \ \ \ \ \ k \in \Omega_L \label{slack1_E}   \\
& -\alpha_kN_k(S_k^{\max}+u_{k,1}^C) \le P_k^C \le \alpha_kN_k(S_k^{\max}+u_{k,2}^C),  \nonumber \\
&\ \ \ \  \ \ \ \ \ \ \ \ \ \ \ \ \ \ \ \ \ \ \ \ \ \ \ \ \ \ \ \ \ \ \ \ \ \ \ \ \ \ \ \ \ \ \ \ \ \ \ k \in \Omega_L^{+} \label{slack1_C}
\end{align}
Constraints (\ref{slack1_E}) and (\ref{slack1_C}) enforce the power flow on the lines that are not connected to zero; however, these two constraints allow thermal violations on the lines in service. The remaining constraints include (\ref{norm_c11})-(\ref{P_B1}), (\ref{Pg_limit1})-(\ref{Pg_fix}).

\section{Case Studies}
\label{case_study}
The proposed planning model is applied to the IEEE 24-bus system and a more practical Polish 2383-bus system. The data for the IEEE 24-bus and the Polish 2383-bus system are included in the MATPOWER software \cite{mybibb:MATPOW}. For all the test systems, each stage is 5 years and all the selected lines and CVSRs are built at the beginning of each stage. The investment cost for the CVSR is assumed to be \$10/kVA \cite{mybibb:Aleks}. Based on the prototype that is going to be installed by Bonneville Power Administration (BPA), the maximum output reactance of the CVSR is allowed to be 20\% of the corresponding line reactance:
\begin{equation}
0 \le x_k^{V} \le 0.2x_k, \ k \in \Omega_{V}  \label{xcvsr_range1}
\end{equation}

\subsection{IEEE 24-Bus System}
\label{IEEE24}
The IEEE 24-bus system has 29 transmission lines, 5 transformers, 32 generators and 21 loads. The thermal limits for all the transmission branches are decreased artificially to introduce congestion. For this test system, we assume only one candidate transmission line per existing line (i.e, excluding transformer upgrades) so the number of candidate transmission lines is 29. In addition, all the existing transmission lines are possible locations to install a CVSR so the number of candidate locations for CVSR is also 29. Excluding one contingency (line 7-8) which splits the system into two parts, complete $N-1$ contingency constraints considering the existing branches are considered. Due to the absence of actual system expansion data, the investment for building new transmission lines is estimated by its length and cost per mile. The cost per mile for different voltage levels can be found in \cite{mybibb:wecc_cost}.

\subsubsection{Single Stage Planning}
We first consider the single stage planning for this test system. The selected lines and CVSRs are committed at the beginning of the stage and the operation cost is evaluated over five years thereafter. The simulation results using integrated model are summarized in Table \ref{24bus_results_inte}. From Table \ref{24bus_results_inte}, it can be seen that the TEP without CVSRs requires building 3 transmission lines. When the CVSR is introduced in the TEP, only 2 transmission lines are needed for the considered stage. The construction of line 14-16 (\$36.47M) is avoided by installing 3 low cost CVSRs (\$13.5M) on line 11-14, 14-16 and 15-21. Thus, the investment cost decreases from \$74.25M to \$51.28M. Although the operating cost of the case with CVSR is \$10M higher than the case without CVSR, the total saving for this five years plan is about \$13 M. The computation time for the case without CVSR is 9.25 s and the time increases to 388.51 s for the case considering CVSR.  

\begin{table}[!htb]
	\centering
	\caption{Single Stage TEP Results Comparison for the IEEE 24-Bus System Using Integrated Model}
	\label{24bus_results_inte}
	\begin{tabular}{|c|| c|c| }
		\hline
		\multirow{2}{*}{}& \multicolumn{2}{c|}{Case}\\
		\cline{2-3}
		&w/o CVSR&w/t CVSR  \\
		\hline
		\multirow{4}{*}{Line}&\multirow{4}{0.7 cm}{14-16 \\ 16-17 \\ 17-18}&\multirow{4}{0.7 cm}{16-17 \\ 17-18} \\
			&&    \\
			&&     \\
			&&      \\
		\hline
	    \multirow{3}{*}{CVSR}&\multirow{3}{*}{-}&\multirow{2}{0.7 cm}{11-14\\ 14-16 \\ 15-21} \\
			&&    \\
			&&    \\
		\hline
		Investment cost (M\$) &74.25&51.28 \\
		\hline
		Operating cost (M\$)  &1168.59&1178.91  \\
		\hline
		Total cost (M\$)  &1242.84&1230.19  \\
		\hline
		Computation Time (s) &9.25&388.51 \\
		\hline
	\end{tabular}
\end{table}

Table \ref{24bus_results_decomposed} shows the TEP results by using the decomposed model. To evaluate the impacts of the decomposition, two cases are simulated: 
\begin{enumerate}
	\item Considering one critical contingency (line 18-21) for the peak and normal load level in the master problem. 
	\item Considering two critical contingencies (line 18-21, 15-21) for the peak and normal load level in the master problem. 
\end{enumerate}       

The critical contingencies are selected based on the circuit loading in the peak load level \cite{mybibb:automatic_contingency}. As observed from Table \ref{24bus_results_decomposed}, the investment plans for the TEP without CVSR are the same for these two cases, which are also identical as the results using integrated model. Nevertheless, the computational time using the decomposed model is only around 1.2 s. The investment plans for the TEP with CVSR are different for the two cases. For the case considering one critical contingency, 1 transmission line and 6 CVSRs are added. The cost in total is \$1234.73M. The case considering two critical contingencies requires to build 2 transmission lines and 3 CVSRs, which are the same planning results as the integrated model. The computational time for the decomposed model considering two critical contingencies is 34.71s. This is 11 times faster than the integrated model. 

\begin{table}[!htb]
	\centering
	\caption{Single Stage TEP Results Comparison for the IEEE 24-Bus System Using Decomposed Model}
	\label{24bus_results_decomposed}
	\begin{tabular}{|c|| c|c|c|c| }
		\hline
		\multirow{2}{*}{}& \multicolumn{2}{c|}{One CC} & \multicolumn{2}{c|}{Two CC}\\
		\cline{2-5}
		&w/o&w/t& w/o&w/t  \\
		&CVSR&CVSR&CVSR&CVSR \\
		\hline
		&14-16&&14-16&\multirow{3}{0.7 cm}{16-17 \\ 17-18}  \\
		Line&16-17&16-17&16-17&  \\
		&17-18&&17-18&  \\
		\hline
		&&11-14&&\multirow{6}{0.7 cm}{11-14 \\ 14-16 \\ 15-21}  \\
		&&14-16&&  \\
		CVSR&-&15-21&-&  \\
		&&17-18&&  \\
		&&17-22&&  \\
		&&21-22&&  \\
		\hline
		Investment&\multirow{2}{*}{74.25}&\multirow{2}{*}{51.28}&\multirow{2}{*}{74.25}&\multirow{2}{*}{51.28}  \\
		cost (M\$)&&&&  \\
		\hline
		Operating  &\multirow{2}{*}{1168.59}&\multirow{2}{*}{1183.45}&\multirow{2}{*}{1168.59}&\multirow{2}{*}{1178.91}  \\
		cost (M\$)&&&&  \\
		\hline
		Total  &\multirow{2}{*}{1242.84}&\multirow{2}{*}{1234.73}&\multirow{2}{*}{1242.84}&\multirow{2}{*}{1230.19}  \\
	    cost (M\$)&&&&  \\
		\hline
		Time (s) &1.15&31.47&1.26&34.71  \\
		\hline
	\end{tabular}
\end{table}

\subsubsection{Multi-stage Planning}
We then consider a two stage planning for this test system. The load growth is estimated to be 25\% in five years and this growth is distributed equally among the load buses. We first evaluate the impacts of $N-1$ contingency constraint on the TEP results. Table \ref{24bus_N_1} summarizes the TEP results with CVSR and without CVSR for the cases that consider and do not consider $N-1$ contingency constraints. The number in the parenthesis indicate the installation year for the new lines and CVSRs. It can be seen that the two cases lead to different network expansion plan. Without CVSR, 3 lines are built for the first stage and no line is needed for the second stage for the case do not consider $N-1$ security constraints. For the case considering $N-1$ security constraints, 2 transmission lines are committed for the first stage and 1 line is added for the second stage. Although the total number of installed transmission lines are the same for the two cases, one long transmission line (15-21) that costs \$69.41M is needed for the case considering $N-1$ security constraint. The construction of this line significantly increases the investment cost for the case considering $N-1$ security constraints. Similar results can also be found in the TEP model with CVSR.          

As observed from Table \ref{24bus_N_1}, for the case considering $N-1$ security constraints, 2 CVSRs on line 11-14 and 14-16 are installed in order to avoid the building of line 14-16. The total savings for this ten years plan is around \$16.63M.    

\begin{table}[!htb]
	\centering
	\caption{Multi-stage TEP Results Comparison for the IEEE-24 Bus System Using Integrated Model}
	\label{24bus_N_1}
	\begin{tabular}{|c|| c|c|c|c| }
		\hline
		\multirow{2}{*}{}& \multicolumn{2}{c|}{Not consider $N-1$} & \multicolumn{2}{c|}{Consider $N-1$}\\
		\cline{2-5}
		&w/o&w/t& w/o&w/t  \\
		&CVSR&CVSR&CVSR&CVSR \\
		\hline
		&14-16 (1)&&14-16 (1)&15-21 (1)  \\
   Lines&16-17 (1)&16-17 (1)&15-21 (1)& 6-10 (6)  \\
		&17-18 (1)&&6-10 (6)&  \\
		\hline
		&&11-14 (1)&&11-14 (1)  \\
	CVSR&-&14-16 (1)&-&14-16 (1)  \\
		&&15-21 (1)&&  \\
		\hline
		Investment&\multirow{2}{*}{74.25}&\multirow{2}{*}{37.78}&\multirow{2}{*}{114.76}&\multirow{2}{*}{87.29}  \\
		cost (M\$)&&&&  \\
		\hline
		Operating  &\multirow{2}{*}{3053.82}&\multirow{2}{*}{3069.82}&\multirow{2}{*}{3049.35}&\multirow{2}{*}{3060.19}  \\
		cost (M\$)&&&&  \\
		\hline
		Total  &\multirow{2}{*}{3128.07}&\multirow{2}{*}{3107.60}&\multirow{2}{*}{3164.11}&\multirow{2}{*}{3147.48}  \\
		cost (M\$)&&&&  \\
		\hline
		Time (s) &0.76&15.07&39.19&790.07  \\
			
		\hline
	\end{tabular}
\end{table}

Table \ref{24bus_muti_decompose} shows the two stage TEP results by using the decomposed model. The same two critical contingencies (line 18-21, 15-21) are considered for the normal and peak load level in stage one and all the load levels in stage two. So the total number of operating states in the master problem is 7 in stage one and 9 in stage two. As observed from Table \ref{24bus_muti_decompose}, the avoidance of building line 14-16 in stage one is achieved by installing 3 CVSRs on line 11-14, 14-16 and 15-21. In addition, the construction of line 18-21 in the second stage is avoided by installing 2 CVSRs on line 18-21 and 21-22. The total saving on the investment is \$44.46M. When comparing the planning results from the integrated model with the results from the decomposed model, one long and expensive transmission line 15-21 (\$69.41M) is installed in stage one in the integrated model. This result arises since forward planning is myopic and does not see the future benefits from the present reinforcement \cite{mybibb:Jabr_TEP,mybibb:TEP_letter}. Still, the difference of the total cost between the decomposed model and the integrated model is \$8.54M for the case considering CVSR, which is only 0.27\% of the planning cost. The computation time of the decomposed model is far less than the integrated model. For the case considering CVSR, the computation time is approximately 18 times faster using the decomposed model.       
\begin{table}[!htb]
	\centering
	\caption{Multi-stage TEP Results Comparison for the IEEE 24-Bus System Using Decomposed Model}
	\label{24bus_muti_decompose}
	\begin{tabular}{|c|| c|c| }
		\hline
		\multirow{2}{*}{}& \multicolumn{2}{c|}{Case}\\
		\cline{2-3}
		&w/o CVSR&w/t CVSR  \\
		\hline
		\multirow{4}{*}{Line}&14-16 (1)& \\
		&16-17 (1)&16-17 (1)    \\
		&17-18 (1)&17-18 (1)     \\
		&6-10 (6)&6-10 (6)      \\
		&18-21 (6)&      \\
		\hline
		&& 11-14 (1)    \\
		&& 14-16 (1)   \\
		CVSR&-&15-21 (1)    \\
		&&18-21 (6)    \\
		&&21-22 (6)    \\
		\hline
		Investment cost (M\$) &111.68&67.22 \\
		\hline
		Operating cost (M\$)  &3059.01&3088.80  \\
		\hline
		Total cost (M\$)  &3170.68&3156.02  \\
		\hline
		Computation Time (s) &2.98&45.13 \\
		\hline
	\end{tabular}
\end{table}

\subsection{Polish 2383-Bus System}
The approach is also applied to a more practical Polish 2383-bus system. The system has 2895 existing branches, 327 generators and 1822 loads. Single stage planning model is used for this case study. Only a few transmission corridors have the potential for the construction with new lines because of the physical or regulatory constraints. It is assumed for this study that the number of candidate lines is 60. In addition, 80 existing transmission lines have been selected as candidate locations to install the CVSR. The selection criterion is the congestion severity of the transmission lines. The line investment cost is estimated by the approach given in Section \ref{IEEE24}. To obtain the contingency list, we first eliminate 643 contingencies that would cause islanding. Then we run an optimal power flow (OPF) for each of the remaining transmission $N-1$ contingencies and select the worst 100 in terms of the operating cost \cite{mybibb:automatic_contingency}. Note that the short term line rating is used for each OPF. Moreover, the worst 6 contingencies are considered for the peak and normal load levels in the master problem. Table \ref{2383bus_decomposed} shows the TEP planning strategy for the case with CVSR and without CVSR by using the decomposed model.

\begin{table}[!htb]
	\centering
	\caption{TEP Results Comparison for the Polish System}
	\label{2383bus_decomposed}
	\begin{tabular}{|c|| c| c| }
		\hline
		\multirow{2}{*}{}& \multicolumn{2}{c|}{Case}\\
		\cline{2-3}
		&w/o CVSR& w/t CVSR  \\
		\hline
		\multirow{8}{*}{Line}&\multirow{8}{2.8 cm}{\centering 437-220, 515-461\\776-539, 1178-834\\994-1289, 1417-1284\\1632-1644, 1693-1632\\1877-1875, 1932-1880\\2328-2165, 2365-2261 \\2348-2379}&\multirow{8}{2.8 cm}{\centering 437-220, 515-461\\776-539, 1178-834\\994-1289, 1417-1284\\1877-1875, 1932-1880\\2328-2165, 2365-2261 \\2348-2379}    \\
		&&    \\
		&&    \\
		&&    \\
		&&    \\
		&&    \\  
		&&    \\
		&&    \\
		\hline
		\multirow{5}{*}{CVSR}&\multirow{5}{*}{\centering -}&\multirow{5}{2.8 cm}{\centering 310-6, 126-127\\613-223, 477-310\\939-1416, 1427-1249\\1693-1632, 1693-1658}      \\
		&&    \\
	    &&    \\
	    &&    \\
	    &&    \\
		\hline
		Investment&\multirow{2}{*}{178.53}&\multirow{2}{*}{182.11}  \\
		cost (M\$)&&  \\
		\hline
		Operating  &\multirow{2}{*}{10474.49}&\multirow{2}{*}{10350.11}  \\
		cost (M\$)&&  \\
		\hline
		Total  &\multirow{2}{*}{10653.02}&\multirow{2}{*}{10532.22}  \\
		cost (M\$)&&  \\
		\hline
		Time (min) &39.23&111.97  \\
		\hline
	\end{tabular}
\end{table}

As observed from Table \ref{2383bus_decomposed}, the TEP without CVSR requires building 13 transmission lines. The total investment cost for this planning strategy is \$178.53M. For the TEP with CVSRs, 11 transmission lines and 8 CVSRs are selected. The investment cost increases by \$3.58M compared to the case without CVSR. Nevertheless, the operating cost decreases significantly by \$124.38M with the inclusion of CVSRs. The saving for this 5 year plan is \$120.9M, which accounts for 1.13\% of the total planning cost. It can also be seen from Table \ref{2383bus_decomposed} that the operating cost takes up a large portion in the total cost for this practical large scale system. The CVSRs are intended to be installed in the appropriate transmission lines to reduce congestion and the operating cost. For the peak load level, the hourly operating cost is \$35988 for the case without CVSRs. The cost is reduced to \$35383 when CVSR is introduced. The computational time when considering CVSR for this practical system is around 1.87 hours, which is acceptable for a planning study \cite{mybibb:tep_wind}.

\subsection{Computational and Optimality Issues}
\label{computational}
The computer used for all simulations has an Intel Core(TM) i5-2400M CPU @ 2.30 GHz with 4.00 GB of RAM. The MILP problem is modeled using the YALMIP \cite{mybibb:YALMIP} toolbox in MATLAB with the CPLEX solver \cite{mybibb:CPLEX} selected to solve the model. 

A heuristic method is used for the decomposed model so the global optimality of the solution is not guaranteed. The impact of the decomposition on the solution can be reduced by including more contingencies in the master problem. This will increase the dimension of the master problem and result in greater computations. So there is a compromise between solution quality and computational time. Table \ref{2383_computational} compares the TEP results for the Polish system considering different number of critical contingencies. As can be observed from the table, TEP considering 6 critical contingencies in the master problem gives better results than TEP considering 3 critical contingencies. Nevertheless, the computational time is higher
for the case considering 6 contingencies. The planner has to balance computational time with solution quality. Note also that each $N-1$ check subproblem takes around 1.2s and is independent. If parallelization techniques are used, the computational time in the subproblem can be significantly reduced and the total time will be largely determined by the master problem. In this case, adding more contingencies in the planning model would have little impact on the computational time as long as the size of the master problem is unchanged, i.e., including the same number of critical contingencies.       

\begin{table}[!htb]
	\centering
	\caption{TEP Results Comparison for the Polish System Using Decomposed Model Considering Different Number of Critical Contingencies}
	\label{2383_computational}
	\begin{tabular}{|c|| c|c|c|c| }
		\hline
		\multirow{2}{*}{}& \multicolumn{2}{c|}{Three CC} & \multicolumn{2}{c|}{Six CC}\\
		\cline{2-5}
		&w/o&w/t& w/o&w/t  \\
		&CVSR&CVSR&CVSR&CVSR \\
		\hline
		No. of Line&15&11&13&11     \\
		\hline
		No. of CVSR&-&5&-&8     \\
		\hline
		Investment&\multirow{2}{*}{191.56}&\multirow{2}{*}{178.99}&\multirow{2}{*}{178.53}&\multirow{2}{*}{182.11}  \\
		cost (M\$)&&&&  \\
		\hline
		Operating  &\multirow{2}{*}{10466.28}&\multirow{2}{*}{10358.02}&\multirow{2}{*}{10474.49}&\multirow{2}{*}{10350.11}  \\
		cost (M\$)&&&&  \\
		\hline
		Total  &\multirow{2}{*}{10657.84}&\multirow{2}{*}{10537.01}&\multirow{2}{*}{10653.02}&\multirow{2}{*}{10532.22}  \\
		cost (M\$)&&&&  \\
		 \hline
		 Master  &\multirow{2}{*}{2.55}&\multirow{2}{*}{48.68}&\multirow{2}{*}{11.34}&\multirow{2}{*}{63.30}  \\
		 problem (min)&&&&  \\
		 \hline
		Total Time (min) &38.59&94.33&39.23&111.97  \\
		\hline
	\end{tabular}
\end{table}

\section{Conclusion}
\label{conclusions}
In this paper, the CVSR is investigated for improving the transmission expansion planning. The CVSR is a FACTS-like device which has the capability of continuously varying the transmission line impedance. Due to its simple and low power rating control circuit, the cost of the CVSR is far less than the cost of a similarly related FACTS device. Thus, a large number of such devices could be installed and have a large impact on the planning process. A security constrained multi-stage TEP model considering CVSR is presented. A reformulation technique is proposed to transform the MINLP model into the MILP model so the model can be efficiently solved by commercial solvers. To relieve the computation burden for a practical large scale system, a decomposition approach is introduced to separate the problem into a planning master problem and security analysis sub-problem. Simulation results on two test systems show that if several CVSRs are appropriately allocated in the system, the building of new transmission lines can be postponed or avoided entirely. Moreover, the CVSRs can change the power flow pattern, which is beneficial in reducing operating cost. Finally, the installation of CVSRs adds flexibility to power system operation and can provide alternative control actions to handle various contingencies.



%




\ifCLASSOPTIONcaptionsoff
  \newpage
\fi



%
\bibliographystyle{IEEEtran}
\bibliography{IEEEabrv,mybibb}

%

\begin{IEEEbiographynophoto}{Xiaohu Zhang}
	(S'12) received the B.S. degree in electrical engineering from Huazhong University of Science and Technology, Wuhan, China, in 2009, the M.S. degree in electrical engineering from Royal Institute of Technology, Stockholm, Sweden, in 2011, and the Ph.D. degree in electrical engineering at The University of Tennessee, Knoxville, in 2017. 
	
	Currently, he works as a power system engineer at GEIRI North America, Santa Clara, CA, USA. His research interests are power system operation, planning and stability analysis.
\end{IEEEbiographynophoto}
\vfill

\begin{IEEEbiographynophoto}{Kevin Tomsovic}
	(F'07) received the B.S. degree in electrical engineering from Michigan Technological University, Houghton, in 1982 and the M.S. and Ph.D. degrees in electrical engineering from the University of Washington, Seattle, in 1984 and 1987, respectively.
	
	Currently, he is the CTI professor of the Department of Electrical Engineering and Computer Science at The University of Tennessee, Knoxville, where he also directs the NSF/DOE ERC CURENT. He was on the faculty of Washington State University from 1992 to 2008. He held the Advanced Technology for Electrical Energy Chair at Kumamoto University, Kumamoto, Japan, from 1999 to 2000 and was an NSF program director in the ECS division from 2004 to 2006.
\end{IEEEbiographynophoto}

\vfill
\begin{IEEEbiographynophoto}{Aleksandar Dimitrovski}
	(SM'09) is the Associate Professor at University of Central Florida – Orlando. Before joining UCF, he had worked in research institutions, power industry, and academia both in USA and Europe. He received his B.Sc. and Ph.D. degrees in electrical engineering with emphasis in power from University Ss. Cyril \& Methodius, Macedonia, and M.Sc. degree in applied computer sciences from University of Zagreb, Croatia. His area of interest is focused on uncertain power systems, their modeling, analysis, protection and control.
\end{IEEEbiographynophoto}







\end{document}